\documentclass[3p,times]{elsarticle}

\usepackage{graphicx}
\usepackage{psfrag}
\usepackage[dvipsnames]{xcolor}

\usepackage{amsmath}
\usepackage{amsthm}
\usepackage{amsfonts}
\usepackage{dsfont}
\usepackage{multirow}
\usepackage{newtxmath}
\usepackage{empheq}

\DeclareMathAlphabet\mathbfcal{OMS}{cmsy}{b}{n}  

\newcommand{\dt}{\Delta t}
\newcommand{\x}{\mathbf{x}}
\renewcommand{\vv}{\mathbf{v}}
\newcommand{\M}{\mathbf{M}}
\newcommand{\R}{\mathbf{R}}
\newcommand{\B}{\mathbf{B}}
\newcommand{\D}{\mathbf{D}}
\newcommand{\C}{\mathbf{C}}
\renewcommand{\R}{\mathbf{R}}
\renewcommand{\S}{\mathbf{S}}

\newcommand{\tav}{\color{black}}
\newcommand{\elo}{\color{black}}
\usepackage{ulem}

\begin{document}
	
	\begin{frontmatter}
		
		\journal{Applied Mathematics and Computation}
		
		

		\title{A high order accurate space-time trajectory reconstruction technique for quantitative particle trafficking analysis}
		
		\author[CIBIO]{Eloina Corradi$^*$}
		\ead{eloina.corradi@unitn.it}
		 \cortext[cor1]{Corresponding author}
        
        \author[unibz]{Maurizio Tavelli}
        
        \author[CIBIO]{Marie-Laure Baudet}
        
        \author[unife]{Walter Boscheri}
        
		\address[CIBIO]{Department of Cellular, Computational and Integrative Biology (CIBIO), University of Trento, Trento 38122, Italy.}
		\address[unibz]{Faculty of Computer Science, Free University of Bozen, Bozen 39100, Italy.}
		\address[unife]{Department of Mathematics and Computer Science, University of Ferrara, Ferrara 44121, Italy.}
%
\begin{abstract}
The study of moving particles (e.g. molecules, virus, vesicles, organelles, or whole cells) is crucial to decipher a plethora of cellular mechanisms within physiological and pathological conditions. Powerful live-imaging approaches enable life scientists to capture particle movements at different scale from cells to single molecules, that are collected in a series of frames. However, although these events can be captured, an accurate quantitative analysis of live-imaging experiments still remains a challenge. Two main approaches are currently used to study particle kinematics: kymographs, which are graphical representation of spatial motion over time, and single particle tracking (SPT) followed by linear linking. Both kymograph and SPT apply a space-time approximation in quantifying particle kinematics, considering the velocity constant either over several frames or between consecutive frames, respectively. Thus, both approaches intrinsically limit the analysis of complex motions with rapid changes in velocity. Therefore, we design, implement and validate a novel reconstruction algorithm aiming at supporting tracking particle trafficking analysis with mathematical foundations. Our method is based on polynomial reconstruction of 4D (3D+time) particle trajectories, enabling to assess particle instantaneous velocity and acceleration, at any time, over the entire trajectory. Here, the new algorithm is compared to state-of-the-art SPT followed by linear linking, demonstrating an increased accuracy in quantifying particle kinematics. Our approach is directly derived from the governing equations of motion, thus it arises from physical principles and, as such, it is a versatile and reliable numerical method for accurate particle kinematics analysis which can be applied to any live-imaging experiment where the space-time coordinates can be retrieved.
\end{abstract}
%
\end{frontmatter}


\section{Introduction} \label{sec.intro}

The use of fluorescent labels and the advance of imaging techniques over the last decades have enabled to capture several cellular and intracellular dynamic events in living samples. The space-time history of the position of moving particles is typically stored in images at different times, which are referred to as frames. Each frame is thus an image that captures particles and their surroundings, at a given time. Along with the improvement of live-imaging techniques, several approaches, algorithms and analysis tools have been implemented to study moving particles, tracking their direction and quantifying their velocity. The current available tools for studying particle trafficking are kymograph and single particle tracking (SPT) followed by linear linking. 

Kymographs are 2D images obtained by stacking the intensity profile along a manually defined path. The result is a graphical representation of spatial motion over time, where the slope of the lines represents the particle velocity along that given path. Kymographs are largely used to investigate particle velocity  \cite{wang2009,daniele2017,lacovich2017,fenn2018,corradi2020}. Although being a widespread technique, kymograph analysis simplifies the complexity of a tracking problem by approximating particle movements along a fixed spatial path, thus two disadvantages intrinsically arise. i) Spatial resolution reduction: by projecting a 3D trajectory onto one single direction and by approximating the kymograph traces with linear segments. ii) Temporal resolution reduction: the dynamics in time is described by tracks with constant behavior (e.g. constant segmental velocity over several acquisitions). All-in-all both space and time resolution are reduced.

Contrary to kymograph, tracking approaches aim at identifying the particle position in each acquisition of a time lapse and at connecting the identified particles over frames, preserving all the spatial information of the dynamics \cite{reid1979,racine2007,serge2008,jaqaman2008,chenouard2014}. Indeed, most of the single-particle tracking (SPT) methods are based on a two-step process: i) detection of particles frame by frame, and ii) particle linking between consecutive frames \cite{serge2008,jaqaman2008}. The final trajectory is obtained with a link of the resulting segments by solving an assignment problem \cite{jonker1987,jaqaman2008}. SPT software are continuously refined in order to increase the performance among different scenarios (e.g. high particle density, low signal-to-noise ratio, motion heterogeneity, blinking signal, merging or splitting events) \cite{cheezum2001,park2010,meijering2012}. Although none of the existing SPT tools tackles all these aspects optimally \cite{chenouard2014}, the selection of one approach or the other enables scientists to study a huge variety of particle dynamics in different biological contexts preserving all the spatial information. Nevertheless, tracking techniques perform a linear linking of identified particles, thus assuming a constant particle velocity between frames. This approximation might represent a limiting step in imaging analysis, especially for particles showing rapid changes in velocity, and it precludes the accurate analysis of particle co-trafficking, i.e. the analysis of a possible correlation of motion among different particles at the same time. 

Therefore we develop and validate a novel algorithm, SHOT-R (Space-time High-Order Trajectory Reconstruction), for high-order particle trajectory reconstruction and kinematic analysis. While SPT approximates the original trajectory by a piecewise first degree polynomial (P1) reconstruction (i.e. connecting dots with linear segments), the SHOT-R algorithm solves particle kinematics equations relying on higher order polynomials. Polynomial reconstruction is a rather general approach in the numerical solution of hyperbolic partial differential equations, mainly focused on advection equations \cite{ADERFSE,IskeKaeser,TumoloBonaventuraRestelli,Shu_ConsSL2011,SLIMEX,boscheri2013semi,MOODorg,ALBI2019_SL}, which is here employed for the first time to the live-imaging context. 
It is unknown how particles move between one acquisition and the next, and any kind of trajectory is an approximation of the real particle path. Nevertheless, here we show through mathematical tests and by applying SHOT-R on a real biological sample, that our method improves the accuracy of the reconstruction and, consequently, the reliability of the kinematic analysis. Indeed, the velocities computed with our algorithm, which accounts for curvilinear trajectories, result to be more accurate than the current SPT strategy. 
Overall, we provide an accurate analysis platform for particle trafficking, which can be applied to any dataset where the space-time coordinates can be retrieved, thus being extremely versatile. Moreover, thanks to the SHOT-R unique feature of a high order continuous trajectory, novel quantitative information can be extracted (i.e. instantaneous velocity and acceleration).

Here, we present an interdisciplinary work that connects state-of-the-art techniques of numerical analysis with cutting-edge live-imaging analysis, that is currently demanding more sophisticated tools to post-process the acquired data. Indeed, the high order accuracy of our method potentially matches the high resolution acquisition techniques that are currently used to maximize the quality of the acquisitions in live-imaging \cite{manzo2015,  lelek2021, parutto2022}. Our method will therefore be presented and validated with the aim of being appreciated from both the mathematical and the biological angle of attack. 
The rest of the paper is organized as follows. In Section \ref{sec.model} we introduce the governing equations of motion, while the numerical scheme is fully detailed in Section \ref{sec.numscheme}. Numerical convergence studies and a large suite of validation tests are then shown in Section \ref{sec.results}. A concluding section finalizes the article where we draw some conclusions and present an outlook to future research. 

\section{Mathematical model} \label{sec.model}
Particle kinematics is described by ordinary differential equations (ODE) governing the motion:
\begin{equation}
	\frac{d \x}{dt} = \vv(\x,t), \qquad \x(t=t_0) = \x_0,
	\label{eqn:ODE}
\end{equation}
where $t$ denotes the time coordinate, $\x=(x,y,z)$ is the position vector of spatial coordinates and $\vv(\x,t)=(v_x,v_y,v_z)$ is the velocity vector with components along each direction in space. The trajectory starts at the initial time $t=t_0$ from the position $\x_0$ (Fig. \ref{fig:Fig1}A) and {\tav then the particle moves in time according to ODE  \eqref{eqn:ODE}}.

The {\tav solution of the} trajectory equation \eqref{eqn:ODE} is a continuous function, whereas experimental data only provide a discrete set of known spatial coordinates $\x_k$ at the corresponding time $t_k$ (Fig. \ref{fig:Fig1}A). Therefore, in order to study particle kinematics, a reconstruction of particle trajectories is needed.

\section{Numerical scheme} \label{sec.numscheme}
\begin{figure}[bhp]
	\centering
	\includegraphics[width=1\linewidth]{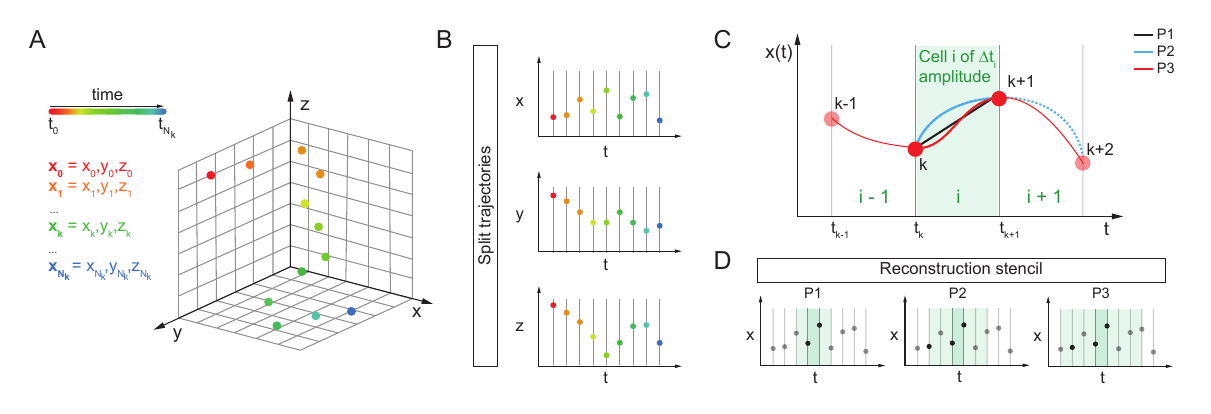}
	\caption{{\bf Space-time high order trajectory reconstruction (SHOT-R) method} {\bf (A)} Schematic of SHOT-R 4D input data: particle motion is described as the evolution of the spatial coordinates $(x,y,z)$ in time (from $t_0$ to $t_{N_k}$). {\bf (B)} Schematic of space-time split trajectory strategy. Each direction $(x,y,z)$ is analyzed separately, thus a multidimensional setting can be reduced to multiple one-dimensional problems. {\bf (C,D)} Schematic of trajectory reconstruction based on polynomials (SHOT-R). An interpolation polynomial $p^N(t)$ of arbitrary degree $N$ is constructed for each cell $i$ (C) by considering the cell itself and a total number $N+1$ surrounding elements (i.e. the reconstruction stencil) (D). Abbreviations: $P1$,$P2$,$P3$, first, second, third degree polynomial.}
	\label{fig:Fig1}
\end{figure}

In the sequel, we detail the derivation and implementation of the SHOT-R method. We assume that a series of frames is available, meaning that we can extract a sequence of points in space and time from a dataset composed of image acquisitions, independently of the nature of the phenomenon under consideration. Obviously, we also assume that the particle motion is always governed by the ODE \eqref{eqn:ODE}. A graphical representation of the reconstruction procedure is shown in Figure \ref{fig:Fig1}.

\subsection{Split trajectories}
In our method, to faithfully describe particle motions, the input data are the space-time coordinates of particles retrieved by the raw movies (i.e. the acquired images of the experiments) using available algorithms for detection and tracking (here TrackMate \cite{tinevez2017trackmate} is employed). 
Once the space-time coordinates $(\x_k,t_k)$ are {\elo retrieved} for each data point $k=1,\ldots,N_K$ (Fig. \ref{fig:Fig1}A), we proceed in a dimension-by-dimension manner by splitting the position vector $\x=(x,y,z)$ along each spatial direction (Fig. \ref{fig:Fig1}B). {\elo This approach speeds up and simplifies the computational process by decreasing the complexity of the problem itself (i.e. from a single 4D problem to multiple 1D ones), and avoiding  the use of a 3D mesh, which is difficult to spatially discretize. At the same time, this strategy keeps the desired temporal and spatial resolution of fully 3D procedures.}

Specifically, the issue related to a multidimensional spatial domain is tackled by considering the trajectory equations	
\begin{equation}
	\label{eqn:ODEsplit}
	\left\{\begin{array}{ll}
		\frac{d x}{dt} = v_x(x,t), & \quad x(t=0) = x_0, \\ 
		\frac{d y}{dt} = v_y(y,t), & \quad y(t=0) = y_0, \\ 
		\frac{d z}{dt} = v_z(z,t), & \quad z(t=0) = z_0,
	\end{array} \right.
\end{equation}
that directly derive from the original definition \eqref{eqn:ODE}.
Let observe that all equations in \eqref{eqn:ODEsplit} are mathematically the same object, therefore {\tav we focus} on the solution of a generic one-dimensional trajectory equation satisfied by a particle along the direction $s$, i.e.
\begin{equation}
	\frac{d s}{dt} = v_s(s,t), \quad s(t=0) = s_0,
\end{equation}
where $s$ can be referred to any spatial direction among $\x$ and the same holds true for the velocity component $v_s$ with $\vv$. This makes our algorithm very general and applicable to a wide range of physical and biological phenomena. We stress that the only information needed as input are the point coordinates $\x_k$ defined along the trajectory that has to be reconstructed.

Let $\mathcal{D}$ represent the spatial dimensionality of the problem under consideration, i.e. $\mathcal{D} \in [1,3]$. Consequently, $\x_k \in \mathds{R}^\mathcal{D}$, that is the spatial coordinates are given by real numbers in dimension $\mathcal{D}$. Along each spatial direction we have the set of space-time coordinates $(s_k,t_k)$ that are extracted from $(\x_k,t_k)$. What carried out for the generic direction $s$ will be performed as many times as the dimensionality of the problem is matched, thus a genuinely multidimensional trajectory defined in dimension $\mathcal{D}$ is reduced to a total number of $\mathcal{D}$ one-dimensional trajectories. The full multidimensional trajectory can then be simply gathered back by collecting all one-dimensional contributions.  

\subsection{Computational mesh}

The time coordinate $t$ is defined in the time interval $T=[t_1;t_{N_K}]$, i.e. $t \in T$, while the space coordinate $s$ is bounded in $\Omega=[s_1;s_{N_K}]$, i.e. $s \in \Omega$. The computational domain $\Omega$ is discretized with a total number of $N_I=N_K-1$ cells with size $\dt_i$ and barycenter time coordinate $t_i$ given by
\begin{equation}
	\label{eqn:mesh}
	\begin{array}{ll}
		\dt_i = t_{i+1/2}-t_{i-1/2}, & \qquad t_i = \frac{t_{i+1/2}+t_{i-1/2}}{2},
	\end{array}
\end{equation}
for $i=1,\ldots,N_I$, where the values at the interfaces are directly available from $t_{i-1/2}=t_k$ and $t_{i+1/2}=t_{k+1}$ (Fig. \ref{fig:Fig1}C). Let observe that index $k$ is used for cell interfaces $i-1/2$ in order to obtain a compatible notation on the computational mesh, thus one has the equivalence $k=i-1/2$. Notice that we do not assume equidistant points, thus the trajectory $s(t)$ can be defined by a sequence of space-time coordinates $(s_k,t_k)$ that are not regularly extracted in time. This does not cause any problem in the algorithm since each cell is assigned its own width $\Delta t_i$. The resulting computational mesh is referred to as staggered mesh, because the known quantities are defined at the interfaces $(k,k+1)=(i-1/2,i+1/2)$ of each cell $i$ whereas the sought unknown reconstruction polynomial $p_i^N(s,t)$ will be discretized at the cell center $s_i$. In other words, known data and unknown reconstruction polynomials adopt staggered definitions on the computational mesh.

\subsection{Reconstruction algorithm} 
The aim of the reconstruction procedure is to obtain a piecewise high order polynomial approximation of the particle trajectory. The result is then a set of polynomials $p_i^N(s,t)$ of degree $N \geq 1$ that are defined within each computational cell $i\in [1,N_I]$ for each spatial direction $s$. The polynomial of arbitrary degree $N$ is written in terms of a normalized Taylor series expanded around the barycenter coordinate $t_i$, that is
\begin{equation}
	p_i^N(s,t) = \sum \limits_{{\tav \ell}=0}^{N} \frac{(t-t_i)^{\tav \ell}}{{\tav \ell}! \, \dt_i^{\tav \ell}} \, \hat{s}_{i{\tav , \ell}} := \phi_\ell(t) \, \hat{s}_{i,\ell}, 
	\label{eqn:polyTaylor}
\end{equation}
with $\dt_i^{\tav \ell}$ being a normalizing factor necessary to avoid ill-conditioned reconstruction matrices and ${\tav \ell}!$ represents as usual the factorial of ${\tav \ell}$. The unknown degrees of freedom are denoted by ${\tav \hat{s}_{i}=\{\hat{s}_{i,\ell}\}_{\ell=0}^N}$ which have to be determined by the reconstruction algorithm. Let $\{\phi_\ell \}_{\ell=0}^{N}$ denote a basis for the polynomial space $\mathds{P}^N$ of degree $N$, the polynomial \eqref{eqn:polyTaylor} can also be written as a linear combination of basis functions $\phi_\ell(t) \in \mathds{P}^N$ belonging to the space of polynomials of degree $N$, i.e. $\mathds{P}^N$. Einstein summation convention is adopted implying summation over repeated indexes, thus the reconstruction polynomial for cell $i$ compactly writes $\sum_\ell \phi_\ell(t) \, \hat{s}_{i,\ell} $ as $\phi_\ell(t) \, \hat{s}_{i,\ell}$.

To retrieve a polynomial of degree $N$, a total number of $N+1$ unknowns $\hat{s}_i$ must be uniquely determined, as follows from the definition \eqref{eqn:polyTaylor}. Depending on the chosen degree $N \geq 1$, the piecewise high order polynomial is reconstructed for each cell $i$ by considering the surrounding known points, that form the so-called reconstruction stencil $\mathcal{N}_k=\{...i-3/2,i-1/2,i+1/2,i+3/2...\}$ (Fig. \ref{fig:Fig1}D). The stencil size depends on the polynomial degree used for the reconstruction (Fig. \ref{fig:Fig1}D). Specifically, the stencil size comprises a total number of cells which is at least equal to the number of degrees of freedom of the chosen polynomial degree. Therefore, the higher the polynomial degree the wider the stencil, i.e. the more information from neighboring cells is required for the reconstruction. The stencil counts a total number of $2(N+1)$ cell interfaces for even and odd degree $N$, so that the stencil is always symmetric with respect to the cell $i$ under consideration. We address with $N$ the polynomial degree, while $\mathcal{O}(N)=N+1$ is the order of accuracy of the method. For instance, a polynomial $P1$ of degree $N=1$ is second order accurate, i.e. $\mathcal{O}(1+1)=2$, because it is exact up to degree $1$ and the approximation errors arise starting from $N=2$, thus the accuracy is of second order. 
The safety factor of $2$ in the total number of cells contained in the stencil {\tav leads to an overdetermined linear system for the degrees of freedom $\hat{s}_i$ that allows for including linear constraints in the reconstruction polynomials.}
 Indeed, we also want to ensure continuity of the reconstructed polynomials across cell interfaces located at $i\pm 1/2$. This would generate a continuous profile for the trajectory, while keeping the reconstruction algorithm local to each cell and therefore computationally efficient.

The reconstruction is designed in such a way that the trajectory coordinates $s_k = s_{i\pm 1/2}$ are exactly matched by the polynomial $p_i^N(s,t)$. For each cell $i$ the reconstruction system is constructed as follows:
\begin{equation}
	\phi_\ell(t_{k}) \,  \hat{s}_{i,\ell} = s_k, \qquad k \in \mathcal{N}_{\tav k},
	\label{eqn:RecSys}
\end{equation}
with the linear constraints
\begin{equation}
	\phi_\ell(t_{i-1/2}) \,  \hat{s}_{i,\ell} = s_{i-1/2}, \qquad 
	\phi_\ell(t_{i+1/2}) \,  \hat{s}_{i,\ell} = s_{i+1/2}.
	\label{eqn:Constr}
\end{equation}
The linear reconstruction system \eqref{eqn:RecSys}-\eqref{eqn:Constr} set up above can be written in {\tav a compact} matrix-vector form as
\begin{equation}
	\left\{ \begin{array}{l} \M \, \S = \B \\ \C \, \S = \D \, \B \end{array} \right.,
	\label{eqn:lsq}
\end{equation}
with the definitions
\begin{eqnarray}
	\B_{[\mathcal{S} \times 1]} &=&\left[ \begin{array}{c}
		\vdots \\ s_{k-1} \\ s_{k} \\ s_{k+1} \\ \vdots 
	\end{array} \right], \qquad 
	\S_{[(N+1) \times 1]} = \left[ \begin{array}{c} \hat{s}_{i,1} \\ \vdots \\ \hat{s}_{i,N+1} \end{array} \right], \qquad 
	\nonumber \\
	\M_{[\mathcal{S} \times (N+1)]} &=& \left[ \begin{array}{ccc}
		\phi_1(t_1) & \ldots & \phi_{N+1}(t_1) \\
		\vdots & \vdots & \vdots \\
		\phi_1(t_k) & \ldots & \phi_{N+1}(t_k) \\
		\vdots & \vdots & \vdots \\
		\phi_1(t_\mathcal{S}) & \ldots & \phi_{N+1}(t_\mathcal{S}) \\
	\end{array} \right].
	\label{eqn:lsq-def}
\end{eqnarray}
The vector of unknown coefficients is $\S$, the right hand side which contains the known position values is $\B$ and the reconstruction matrix is given by $\M$. The constrained matrix $\C$ in \eqref{eqn:lsq} is a sub-partition of $\M$ and contains those rows which correspond to the stencil elements $i \pm 1/2$, while the vector $\D$ is zero everywhere apart from the two entries related to the constraints \eqref{eqn:Constr}, where it is set to $1$. The constrained least squares (CLSQ) method is adopted for solving system \eqref{eqn:lsq}, therefore a functional $g(\S)$ is defined and minimized by requiring its derivatives to vanish, thus
\begin{eqnarray}
	g(\S) &=& \left( \M { \mathbf{S}} - \B \right)^T \left( \M { \mathbf{S}} - \B \right) - \boldsymbol{\mu}^T \left( \mathbf{C} { \mathbf{S}} - \mathbf{D} \B \right), \\
	\frac{\partial g}{\partial \S} &=& 2 \mathbf{M} \mathbf{M} \S - 2 \mathbf{M} \B - \mathbf{C} \boldsymbol{\mu} = 0, \\
	\frac{\partial g}{\partial \boldsymbol{\mu}} &=& -\left( \mathbf{C} { \mathbf{S}} - \mathbf{D} \B \right) = 0,
\end{eqnarray}
where $\boldsymbol{\mu}$ is the vector of Lagrange multipliers to enforce the linear constraints. The associated enlarged linear system of normal equations then reads
\begin{equation}
	\left( \begin{array}{cc} 2 \mathbf{M} \mathbf{M} & -\mathbf{C} \\ \mathbf{C} & \mathbf{0} \end{array} \right) \left( \begin{array}{c} \S \\ \boldsymbol{\mu} \end{array} \right) = \left( \begin{array}{c} 2 \mathbf{M} \\ \mathbf{D} \end{array} \right) \B,
	\label{eqn:clsq-1}
\end{equation}
whose solution can be simply obtained as
\begin{equation}
	\left( \begin{array}{c} \S \\ \boldsymbol{\mu} \end{array} \right) = \left( \begin{array}{cc} 2 \mathbf{M} \mathbf{M} & -\mathbf{C} \\ \mathbf{C} & \mathbf{0} \end{array} \right)^{-1} \left( \begin{array}{c} 2 \mathbf{M} \\ \mathbf{D} \end{array} \right) \, \B := \R^L \, \B. 
	\label{eqn:clsq}
\end{equation}
We are only interested in the solution of the expansion coefficients $\S$, therefore the final reconstruction matrix $\R$ for the trajectory in the $s$ dimension is a subset of the CLSQ matrix $\R^L$, that is $\R:=\R^L_{[(N+1) \times \mathcal{N}_i]}$. Once the reconstruction matrices have been computed for all cells, the reconstruction algorithm is efficiently carried out by one matrix-vector product, namely
\begin{equation}
	\S = \R^L_{[(N+1) \times \mathcal{N}_i]} \, \B,
	\label{eqn:recS}
\end{equation}
which permits to obtain the sought degrees of freedom $\hat{s}_{i,\ell}$ and thus to uniquely define the reconstruction polynomial $p_i^N(s,t)$ according to \eqref{eqn:polyTaylor} within each cell $i$.

\subsection{Nonlinear limiting with CWENO strategy}
\label{cweno}
In the presence of discontinuities in the trajectory, which might be due either to the real particle motion or to the casualty of the acquisition process, a limiting technique must be introduced in the reconstruction procedure. The limiter is needed to avoid the generation of spurious oscillations that arise when computing a high order interpolation of the trajectory. Since high order polynomials are not monotone, the reconstruction polynomial $p_i^N(s,t)$ could describe nonphysical oscillations that we want to eliminate. To that aim, we rely on the Central Weighted Essentially Non-Oscillatory (CWENO) strategy firstly proposed in \cite{LPR:99} for hyperbolic systems of conservation laws. 

Differently from \cite{LPR:99}, here we apply the CWENO strategy in a pointwise manner, hence resorting to a kind of finite difference WENO scheme \cite{ShuFD2003}. In the CWENO framework the polynomial $p_i^N(s,t)$ given by \eqref{eqn:polyTaylor} is called \textit{optimal polynomial}, because among all the possible polynomials of degree $N$, it is close in the least-square sense to the point values $s_{k}$ in the stencil $\mathcal{N}_k$, while satisfying the linear constraints \eqref{eqn:Constr}. This reconstruction is linear in the sense of Godunov \cite{Godunov1959}, thus it must be limited and stabilized via a nonlinear scheme. The nonlinearity of the reconstruction lies in the weights used to combine different reconstruction polynomials. According to \cite{LPR:99,ADER_CWENO}, the central polynomial $p_{i,0}^{N}(s,t)$ is evaluated as the difference between the optimal polynomial $p_i^N(s,t)$ and the linear combination of the two one-sided polynomials $p_{i,1}^{1}(s,t)$ (left) and $p_{i,2}^{1}(s,t)$ (right) of degree $N=1$ \cite{CPSV:2018}, that is
\begin{equation}
	\label{CWENO:P0}
	p_{i,0}^{N}(s,t)= \frac{1}{\lambda_0}\left(p_i^{N}(s,t) - \lambda_{1} p_{i,1}^{1}(s,t) - \lambda_{2} p_{i,2}^{1}(s,t) \right),
\end{equation}
where $\lambda_{0},\lambda_{1}, \lambda_2$ are positive coefficients such that
\begin{equation}
	\lambda_{0} + \lambda_{1} + \lambda_2=1.
	\label{eqn.sumCWENO}
\end{equation}
The one-sided second order polynomials are simply computed as the linear polynomial connecting the point value $s_k$ with the direct left or right neighbor point, namely $s_{k-1}$ or $s_{k+1}$, respectively:
\begin{equation}
	p_{i,1}^{1}(s,t) = s_{k} + \frac{s_k-s_{k-1}}{t_k - t_{k-1}} \cdot t, \qquad p_{i,2}^{1}(s,t) = s_{k} + \frac{s_{k+1}-s_{k}}{t_{k+1}-t_k} \cdot t.
	\label{eqn.P1}
\end{equation}
Let us notice that the above relations can still be represented in terms of the expansion given by \eqref{eqn:polyTaylor} by setting $N=1$. These lower polynomials are nothing but the equivalent of the linear linking used in SPT methods. Here, the second order reconstructions \eqref{eqn.P1} are instead used to limit the central high order reconstruction polynomial $p_{i,0}^{N}(s,t)$. The linear weights in \eqref{eqn.sumCWENO} are a normalization which sums up to unity and we set $\lambda_{0}=200/\lambda_{sum}$ for the central polynomial, while we take $\lambda_{1}=\lambda_{2}=1/\lambda_{sum}$ for the left and right one-sided polynomials, hence $\lambda_{sum}=202$. The nonlinear data-dependent hybridization among the three polynomials obtained for each stencil is given by
\begin{equation}
	\label{CWENOx}
     \tilde{p}_i^{N}(s,t) = \omega_0 \, p_{i,0}^{N}(s,t) + \omega_{1} p_{i,1}^{1,\ell}(s,t) + \omega_{2} p_{i,2}^{1,\ell}(s,t),  
\end{equation}
where the nonlinear weights $\{\omega_m\}_{m=0}^2$ are given by 
\begin{equation}
	\label{eqn.weights}
	\omega_m = \frac{\tilde{\omega}_m}{\sum \limits_{m=0}^{2} \tilde{\omega}_m}, 
	\qquad \textnormal{ with } \qquad 
	\tilde{\omega}_m = \frac{\lambda_m}{\left(\sigma_m + \epsilon \right)^r}. 
\end{equation} 
The parameter $\epsilon=10^{-14}$ avoids division by zero and the exponent $r=4$ is chosen according to \cite{ADER_CWENO}. The oscillation indicators $\sigma_m$ are computed according to \cite{DumbserEnauxToro} as
\begin{equation}
	\sigma_m = \sum_{\alpha=1}^{N} \int_{t_k}^{t_{k+1}} \frac{\partial^\alpha \phi_\ell(t)}{\partial t^\alpha} \, \hat{s}_{m,\ell} \cdot \frac{\partial^\alpha \phi_r(t)}{\partial t^\alpha} \, \hat{s}_{m,r} d t. 
\end{equation} 
Eventually, the unlimited reconstruction polynomial \eqref{eqn:polyTaylor} is overwritten as follows:
\begin{equation}
	{p}_i^{N}(s,t) = \tilde{p}_i^{N}(s,t),
	\label{eqn.CWENOpoly}
\end{equation}
which results in a stable and limited high order approximation.

\subsection{Computation of the length of curvilinear trajectory}
\label{isoP}
The evaluation of the length $L$ of the particle trajectory is needed in order to compute average quantities over the entire path. The length is defined as the integral along the line drawn by the particle trajectory. The integral can be rewritten as the sum of the integrals over all computational cells $N_I$, that is
\begin{equation}
	L = \int \limits_{\x_1}^{\x_{N_K}} \, dl = \sum_{i=1}^{N_I} \int \limits_{\x_{i-1/2}}^{\x_{i+1/2}} \, dl.
	\label{eqn:IntL}
\end{equation}
The coordinates $(\x_1,\x_{N_K})$ denote the starting and the ending point of the trajectory in the space of dimension $\mathcal{D} \in [1;3]$. To perform the integration \eqref{eqn:IntL} at high order of accuracy, it is not sufficient to split the trajectory into $N_I$ segments according to the number of cells and to sum up all the linear contributions. Indeed, this would lead to a second order accurate evaluation of $L$, since piecewise linear polynomials $P1$ are used { and \eqref{eqn:IntL} reduces to the classical trapezoidal rule}.

To obtain a consistent high order computation of the trajectory length, isoparametric discretizations must be adopted along the generic spatial direction $s$. To this aim, let us first define a reference coordinate system $\xi \in [0;1]$ which is used to map a generic cell $i$, hence
\begin{equation}
	s = s_{i-1/2} + \xi \left( s_{i+1/2} - s_{i-1/2} \right),
	\label{eqn:map1d}
\end{equation}
so that for $\xi=0$ the left interface coordinate $s_{i-1/2}$ of the cell is retrieved, while for $\xi=1$ the right interface $s_{i+1/2}$ is obtained. In the reference system we then define a set of nodal basis functions $\theta_m(\xi)$ which are chosen to be the Lagrange polynomials passing through the set of nodes $\hat{\xi}_m=m/N$ with $0 \leq m \leq N$. Furthermore, the
standard interpolation condition $\theta_q(\hat{\xi}_m) = \delta_{qm}$ holds for the basis, with $\delta_{qm}$ being the usual Kronecker symbol. Up to degree $N=3$ the nodal basis functions $\theta_m(\xi)$ with associated nodes $\hat{\xi}_m$ explicitly write as follows: 
\begin{itemize}
	\item $N=1$
	\begin{equation}
		\begin{array}{ll}
			\theta_1 = 1 - \xi &
			\theta_2 = \xi \\
			& \\
			\hat{\xi}_1 = 0 &
			\hat{\xi}_2 = 1  
		\end{array}
		\label{eqn:N1} 
	\end{equation}
	\item $N=2$
	\begin{equation}
		\begin{array}{lll}
			\theta_1 = 1 - 3 \xi + 2 \xi^2 &
			\theta_2 = 4 \xi - 4 \xi^2 &
			\theta_3 = -\xi + 2 \xi^2 \\
			& & \\
			\hat{\xi}_1 = 0 &
			\hat{\xi}_2 = 1/2 &  
			\hat{\xi}_3 = 1 \\
		\end{array}
		\label{eqn:N2}  
	\end{equation}
	\item $N=3$
	\begin{equation}
		\begin{array}{ll}
			\theta_1 = 1 - \frac{11}{2}  \xi + 9  \xi^2 - \frac{9}{2}  \xi^3 &
			\theta_2 = 9  \xi - \frac{45}{2}  \xi^2 + \frac{27}{2}  \xi^3 \\
			\theta_3 = - \frac{9}{2}  \xi + 18  \xi^2 - \frac{27}{2}  \xi^3 &
			\theta_4 = \xi - \frac{9}{2}  \xi^2 + \frac{9}{2}  \xi^3 \\
			& \\
			\hat{\xi}_1 = 0 &
			\hat{\xi}_2 = 1/3 \\  
			\hat{\xi}_3 = 2/3 &
			\hat{\xi}_4 = 1 
		\end{array}  
		\label{eqn:N3}
	\end{equation}
\end{itemize}
The integral over cell $i$ in \eqref{eqn:IntL} can be performed in the reference system at the aid of a change of variables $\x \to \xi$, hence obtaining
\begin{equation}
	\int \limits_{\x_{i-1/2}}^{\x_{i+1/2}} \, dl = \int \limits_{0}^{1} |J_\xi| \, d\xi,
	\label{eqn:IntL2}
\end{equation}
with the Jacobian determinant given by
\begin{equation}
	|J_\xi| = \sqrt{\left(\frac{\partial x }{\partial \xi}\right)^2 + \left(\frac{\partial y }{\partial \xi}\right)^2 + \left(\frac{\partial z }{\partial \xi}\right)^2}.
	\label{eqn:detJ}
\end{equation}
The computation of the derivatives in \eqref{eqn:detJ} is done as usual by exploiting the expansion of the basis functions $\theta(\xi)$, thus for a generic spatial coordinate $s$ it holds 
\begin{equation}
	\frac{\partial s }{\partial \xi} = \frac{\partial \theta_m(\xi) }{\partial \xi} \hat{s}_m, \label{eqn:dsdxi}
\end{equation}
where the degrees of freedom $\hat{s}_m$ represent the value of the coordinate $s$ corresponding to the reference node $\hat{\xi}_m$. To be precise, the values $\hat{s}_m$ are computed according to \eqref{eqn:evalX}, relying on the SHOT-R reconstruction polynomial. Finally, the integration in \eqref{eqn:IntL2} is numerically approximated by Gaussian quadrature formulae of suitable order of accuracy \cite{stroud}.

\subsection{Computation of position, velocity and acceleration}

The reconstruction polynomials $p_i^N(s,t)$ provide a piecewise continuous representation of the trajectory within each cell $i$ along the generic spatial direction $s$ as a function of time $t$. Thus, the kinematic description (position $s$, velocity $v_s$ and acceleration $a_s$) of the particle at any given time $\tilde{t}$ can be easily computed in two steps:

\begin{itemize} 
	\item find the computational cell containing $\tilde{t}$, i.e. {$i\in[1,N_I]$ such that} $t_{i-1/2} \leq \tilde{t}\leq t_{i+1/2}$, so that the associated reconstruction polynomial $p_i^N(s,t)$ can be identified; \\
	\item extract the kinematic quantities 
	\begin{eqnarray}
		s(\tilde{t}) &= p_i^N(s,\tilde{t}) &=  \sum \limits_{q=0}^{N} \frac{(\tilde{t}-t_i)^q}{q! \, \dt_i^q} \, \hat{s}_i,  \label{eqn:evalX} \\
		v_s(\tilde{t}) &= \frac{d \, p_i^N(s,\tilde{t})}{dt} &=  \sum \limits_{q={\tav 1}}^{N} \frac{q(\tilde{t}-t_i)^{q-1}}{q! \, \dt_i^q} \, \hat{s}_i, \label{eqn:evalV} \\
		a_s(\tilde{t}) &= \frac{d^2 \, p_i^N(s,\tilde{t})}{dt^2} &=  \sum \limits_{q={\tav 2}}^{N} \frac{q(q-1)(\tilde{t}-t_i)^{q-2}}{q! \, \dt_i^q} \, \hat{s}_i. \label{eqn:evalA}
	\end{eqnarray}
\end{itemize} 
The above definitions are referred to as instantaneous quantities. We also need to give explicit formulation for the global average speed $v_L$, the displacement velocity $v_D$ and the mean of instantaneous velocities within each frame $\bar{v}_M$. These definitions are mathematically described as follows:
\begin{eqnarray}
	v_L &=& \frac{L}{t_{N_K}-t_1}, \qquad L=\int_{\x_1}^{\x_{N_K}} dl,  \label{eqn:vL}\\
	v_D &=& \frac{s_{N_K}-s_{1}}{t_{N_K}-t_1}  \label{eqn:vD}, \\
	v_M &=& \frac{1}{N_I} \sum_{i=1}^{N_I} v_{i}, \qquad v_i = \frac{s_{i+1}-s_i}{t_{i+1}-t_i}.  \label{eqn:vM}
\end{eqnarray}
with $s$ denoting as usual a generic spatial direction among $\x=(x,y,z)$. The length of the trajectory $L$ in \eqref{eqn:vL} is computed as the integral along the curvilinear path from $\x_1$ to $\x_{N_K}$ according to \eqref{eqn:IntL}. Finally, the instantaneous velocity at a generic space-time coordinate $(s,\tilde{t})$ can be simply evaluated relying on the velocity reconstruction polynomial, thus according to \eqref{eqn:evalV}. 

In order to fully exploit the piecewise continuous information provided by the high order reconstruction, the evaluation of all instantaneous kinematic quantities (position, velocity and acceleration) is extracted on a very fine computational mesh given by $N+1$ Gauss points within each computational cell, hence leading to exact integration up to order $2(N+1)-1$ (see \cite{stroud}) within each cell.

\section{Numerical results} \label{sec.results}

The test suite is chosen such that we rigorously validate the accuracy of our method from the mathematical viewpoint, but we also propose some test cases to demonstrate the applicability of the algorithm to biological datasets. 

All data are analyzed with Prism (GraphPad 6 or 7) {\elo and for} all statistical tests the significance level is $\alpha = 0.05$. The exact number of replicates, the type of statistical tests {\elo and p-values} are reported in Figure legends.

\subsection{Numerical convergence studies} 
The order of accuracy $p > 0$ of a numerical method is also defined as the largest value of $p$ for which the following inequality holds \cite{LeVequeODE}:
\begin{equation}
	\varepsilon(\dt) \leq C \, \dt^p, 
	\label{eqn.Accuracy}
\end{equation}
where $\varepsilon(\dt)$ is the error of the numerical solution, which depends on the mesh size $\dt$. The constant $C{\tav >0}$ {\tav depends on the problem but} is independent of $\dt$, this is also written thus
\begin{equation}
	\varepsilon = \mathcal{O}(\Delta t^p), \qquad \textnormal{ as } \quad \dt \to 0.
\end{equation}
In other words, the numerical error must decrease with convergence rate of $\mathcal{O}(\Delta t^p)$ as the mesh size decreases. Given two regular meshes with mesh sizes $\dt_1$ and $\dt_2$, respectively, and the same numerical method, we have a fixed relationship between the meshes which writes
\begin{eqnarray}
	\varepsilon(\dt_1) &=& C\dt_1^p, \label{eqn.errDT1} \\
	\varepsilon(\dt_2) &=& C\dt_2^p. \label{eqn.errDT2}
\end{eqnarray}
Dividing \eqref{eqn.errDT2} by \eqref{eqn.errDT1} we obtain
\begin{equation}
	\left( \frac{\dt_2}{\dt_1}\right)^p = \frac{\varepsilon(\dt_2)}{\varepsilon(\dt_1)},
\end{equation}
which leads to
\begin{equation}
	p = \frac{\ln \left(\frac{\varepsilon(\dt_2)}{\varepsilon(\dt_1)}\right)}{\ln \left(\frac{\dt_2}{\dt_1}\right)},
\end{equation}
that is the empirical order of accuracy. It must be $p=N+1$, meaning that a numerical scheme with reconstructions of degree $N$ introduces numerical errors starting from degree $N+1$, that is a polynomial of degree $N$ can be exactly reconstructed. In practice one performs numerical experiments for a sequence of refined meshes $\{\dt_1,\dt_2, \ldots ,\dt_K\}$ with usually a fixed ratio between successive meshes of sizes $\dt_k$ and $\dt_{k+1}$. For example $\dt_{k+1}=\frac{\dt_{k}}{2}$. A sequence of numbers $\{p_2,p_3,\ldots,p_K\}$ is then retrieved which is given by
\begin{equation}
	p_k = \frac{\ln \left( \frac{\varepsilon(\dt_k)}{\varepsilon(\dt_{k-1})}\right)}{\ln \left( \frac{\dt_k}{\dt_{k-1}}\right)}, \qquad k=2,\ldots,K.
\end{equation}	

To carry out the convergence studies, we consider the time interval $T=[-1;1]$ and we prescribe the following analytical trajectory in 3D:
\begin{equation}
	\x=(x,y,z) = \left( \sin(\pi x) \cos(2\pi x), \, 
	3\cos(2\pi y) - 2\sin(\pi y), \, -6\sin(\pi z)+2\cos(3\pi z) \right).
	\label{eqn:IniConv}
\end{equation} 
The new SHOT-R method is then employed to reconstruct the trajectory along each direction and the errors $\varepsilon$ are measured for each spatial direction $(x,y,z)$ in $L_1$, $L_2$ and $L_\infty$ norms according to \eqref{eqn:Lnorms}.
The exact solution is explicitly computed with \eqref{eqn:IniConv}, while the reconstructed position is evaluated according to \eqref{eqn:evalX}. A sequence of four successively refined meshes is used with $\dt=\{100,200,400,800\}$ and the results are reported in Table \ref{tab:convTable}. The errors are measured in the error norms described in \ref{app.err_norms} (Fig. \ref{fig:Fig2}A), and the results numerically confirm that the SHOT-R algorithm is actually high order accurate ($N+1$) for an arbitrary polynomial degree $N$.

\begin{center} 
	\begin{table}[!htbp] 
			\begin{footnotesize}
				\setlength{\tabcolsep}{0.3pt}
				\caption{{\bf Convergence table} Numerical convergence rates obtained with high order accurate reconstruction of the velocity field given in Section B of the article. Errors in $L_1$, $L_2$ and $L_\infty$ norms are reported for all position components $(x,y,z)$ with associated order of accuracy $\mathcal{O}$ on a sequence of refined mesh with characteristic size $\dt$.}
				\begin{tabular}{|c|cccccc|cccccc|cccccc|}
					\multicolumn{19}{c}{}  \\
					\multicolumn{19}{c}{$\mathbf{N=1}$} \\
					\hline
					\hline
					\multicolumn{19}{c}{} \\
					\multicolumn{7}{c}{horizontal position $x$} & \multicolumn{6}{c}{horizontal position $y$} & \multicolumn{6}{c}{horizontal position $z$} \\
					\hline
					$\dt$ & $\varepsilon_{L_1}$ & $\mathcal{O}(L_1)$ & $\varepsilon_{L_2}$ & $\mathcal{O}(L_2)$ & $\varepsilon_{L_\infty}$ & $\mathcal{O}(L_\infty)$ & $\varepsilon_{L_1}$ & $\mathcal{O}(L_1)$ & $\varepsilon_{L_2}$ & $\mathcal{O}(L_2)$ & $\varepsilon_{L_\infty}$ & $\mathcal{O}(L_\infty)$ &$\varepsilon_{L_1}$ & $\mathcal{O}(L_1)$ & $\varepsilon_{L_2}$ & $\mathcal{O}(L_2)$ & $\varepsilon_{L_\infty}$ & $\mathcal{O}(L_\infty)$ \\ 
					\hline
					1.00e-2 & 1.89e-3 & -     & 1.49e-3 & -     & 1.64e-3 & -    & 5.06e-3 &  -    & 4.00e-3 & -     & 4.60e-3 & -    & 7.74e-3 & -     & 6.24e-3 & -     & 7.62e-3 & -    \\
					5.00e-3 & 4.73e-4 & 2.00  & 3.72e-4 & 2.00  & 4.11e-4 & 2.00 & 1.27e-3 & 2.00  & 1.00e-3 & 2.00  & 1.15e-3 & 2.00 & 1.94e-3 & 2.00  & 1.56e-3 & 2.00  & 1.91e-3 & 2.00 \\
					2.50e-3 & 1.18e-4 & 2.00  & 9.31e-5 & 2.00  & 1.03e-4 & 2.00 & 3.16e-4 & 2.00  & 2.50e-4 & 2.00  & 2.88e-4 & 2.00 & 4.84e-4 & 2.00  & 3.90e-4 & 2.00  & 4.77e-4 & 2.00 \\
					1.25e-3 & 2.95e-5 & 2.00  & 2.33e-5 & 2.00  & 2.57e-5 & 2.00 & 7.91e-5 & 2.00  & 6.25e-5 & 2.00  & 7.20e-5 & 2.00 & 1.21e-4 & 2.00  & 9.75e-5 & 2.00  & 1.19e-4 & 2.00 \\
					\hline
					\multicolumn{19}{c}{}  \\
					
					\multicolumn{19}{c}{$\mathbf{N=2}$} \\
					\hline
					\hline
					\multicolumn{19}{c}{} \\
					\multicolumn{7}{c}{horizontal position $x$} & \multicolumn{6}{c}{horizontal position $y$} & \multicolumn{6}{c}{horizontal position $z$} \\
					\hline
					$\dt$ & $\varepsilon_{L_1}$ & $\mathcal{O}(L_1)$ & $\varepsilon_{L_2}$ & $\mathcal{O}(L_2)$ & $\varepsilon_{L_\infty}$ & $\mathcal{O}(L_\infty)$ & $\varepsilon_{L_1}$ & $\mathcal{O}(L_1)$ & $\varepsilon_{L_2}$ & $\mathcal{O}(L_2)$ & $\varepsilon_{L_\infty}$ & $\mathcal{O}(L_\infty)$ &$\varepsilon_{L_1}$ & $\mathcal{O}(L_1)$ & $\varepsilon_{L_2}$ & $\mathcal{O}(L_2)$ & $\varepsilon_{L_\infty}$ & $\mathcal{O}(L_\infty)$ \\
					\hline
					1.00e-2 & 2.80e-4 & -     & 2.49e-4 & -     & 6.16e-4 & -    & 4.86e-4 & -     & 4.14e-4 & -     & 5.91e-4 & -    & 1.11e-3 & -     & 9.49e-4 & -     & 1.38e-3 & -    \\
					5.00e-3 & 3.43e-5 & 3.03  & 3.02e-5 & 3.04  & 8.16e-5 & 2.92 & 5.96e-5 & 3.03  & 5.13e-5 & 3.01  & 7.39e-5 & 3.00 & 1.35e-4 & 3.04  & 1.16e-4 & 3.03  & 1.72e-4 & 3.00 \\
					2.50e-3 & 4.23e-6 & 3.02  & 3.69e-6 & 3.03  & 1.03e-5 & 2.98 & 7.42e-6 & 3.01  & 6.41e-6 & 3.00  & 9.24e-6 & 3.00 & 1.68e-5 & 3.01  & 1.45e-5 & 3.01  & 2.15e-5 & 3.00 \\
					1.25e-3 & 5.25e-7 & 3.01  & 4.55e-7 & 3.02  & 1.30e-6 & 2.99 & 9.27e-7 & 3.00  & 8.01e-7 & 3.00  & 1.15e-6 & 3.00 & 2.09e-6 & 3.00  & 1.81e-6 & 3.00  & 2.69e-6 & 3.00 \\
					\hline
					\multicolumn{19}{c}{}  \\
					
					\multicolumn{19}{c}{$\mathbf{N=3}$} \\
					\hline
					\hline
					\multicolumn{19}{c}{} \\
					\multicolumn{7}{c}{horizontal position $x$} & \multicolumn{6}{c}{horizontal position $y$} & \multicolumn{6}{c}{horizontal position $z$} \\
					\hline
					$\dt$ & $\varepsilon_{L_1}$ & $\mathcal{O}(L_1)$ & $\varepsilon_{L_2}$ & $\mathcal{O}(L_2)$ & $\varepsilon_{L_\infty}$ & $\mathcal{O}(L_\infty)$ & $\varepsilon_{L_1}$ & $\mathcal{O}(L_1)$ & $\varepsilon_{L_2}$ & $\mathcal{O}(L_2)$ & $\varepsilon_{L_\infty}$ & $\mathcal{O}(L_\infty)$ &$\varepsilon_{L_1}$ & $\mathcal{O}(L_1)$ & $\varepsilon_{L_2}$ & $\mathcal{O}(L_2)$ & $\varepsilon_{L_\infty}$ & $\mathcal{O}(L_\infty)$ \\
					\hline
					1.00e-2 & 7.66e-5 & -     & 6.60e-5 & -     & 1.04e-4 & -    & 9.05e-5 &  -    & 8.24e-5 & -     & 2.27e-4 & -    & 2.99e-4 & -     & 2.69e-4 & -     & 6.88e-4 & -    \\
					5.00e-3 & 4.68e-6 & 4.03  & 4.01e-6 & 4.04  & 4.94e-6 & 4.39 & 5.61e-6 & 4.01  & 5.00e-6 & 4.04  & 1.50e-5 & 3.92 & 1.88e-5 & 3.99  & 1.68e-5 & 4.00  & 4.94e-5 & 3.80 \\
					2.50e-3 & 2.90e-7 & 4.01  & 2.50e-7 & 4.00  & 3.10e-7 & 3.99 & 3.47e-7 & 4.01  & 3.05e-7 & 4.04  & 9.51e-7 & 3.98 & 1.17e-6 & 4.01  & 1.03e-6 & 4.03  & 3.19e-6 & 3.95 \\
					1.25e-3 & 1.81e-8 & 4.00  & 1.56e-8 & 4.00  & 1.94e-8 & 4.00 & 2.16e-8 & 4.01  & 1.88e-8 & 4.02  & 5.96e-8 & 4.00 & 7.29e-8 & 4.01  & 6.34e-8 & 4.02  & 2.01e-7 & 3.99 \\
					\hline
					\multicolumn{19}{c}{}  \\
					
					\multicolumn{19}{c}{$\mathbf{N=4}$} \\
					\hline
					\hline
					\multicolumn{19}{c}{} \\
					\multicolumn{7}{c}{horizontal position $x$} & \multicolumn{6}{c}{horizontal position $y$} & \multicolumn{6}{c}{horizontal position $z$} \\
					\hline
					$\dt$ & $\varepsilon_{L_1}$ & $\mathcal{O}(L_1)$ & $\varepsilon_{L_2}$ & $\mathcal{O}(L_2)$ & $\varepsilon_{L_\infty}$ & $\mathcal{O}(L_\infty)$ & $\varepsilon_{L_1}$ & $\mathcal{O}(L_1)$ & $\varepsilon_{L_2}$ & $\mathcal{O}(L_2)$ & $\varepsilon_{L_\infty}$ & $\mathcal{O}(L_\infty)$ &$\varepsilon_{L_1}$ & $\mathcal{O}(L_1)$ & $\varepsilon_{L_2}$ & $\mathcal{O}(L_2)$ & $\varepsilon_{L_\infty}$ & $\mathcal{O}(L_\infty)$ \\
					\hline
					1.00e-2 & 1.18e-5 & -     & 1.17e-5 &  -    & 4.77e-5 & -    & 9.72e-6 & -     & 8.65e-6 & -     & 2.74e-5 & -    & 5.11e-5 & -     & 4.89e-5 & -     & 1.88e-4 & -    \\
					5.00e-3 & 3.68e-7 & 5.01  & 3.67e-7 & 4.99  & 1.98e-6 & 4.59 & 2.76e-7 & 5.14  & 2.36e-7 & 5.20  & 4.66e-7 & 5.88 & 1.42e-6 & 5.17  & 1.23e-6 & 5.32  & 3.31e-6 & 5.83 \\
					2.50e-3 & 1.11e-8 & 5.05  & 1.06e-8 & 5.12  & 6.60e-8 & 4.91 & 8.40e-9 & 5.04  & 7.21e-9 & 5.03  & 9.98e-9 & 5.55 & 4.27e-8 & 5.06  & 3.66e-8 & 5.07  & 5.48e-8 & 5.92 \\
					1.25e-3 & 3.38e-10 & 5.04  & 3.09e-10 & 5.10  & 2.09e-9 & 4.98 & 2.61e-10 & 5.01  & 2.25e-10 & 5.00  & 3.12e-10 & 5.00 & 1.32e-9 & 5.01  & 1.14e-9 & 5.01  & 1.57e-9 & 5.12 \\
					\hline
					\multicolumn{19}{c}{}  \\
					
					\multicolumn{19}{c}{$\mathbf{N=5}$} \\
					\hline
					\hline
					\multicolumn{19}{c}{} \\
					\multicolumn{7}{c}{horizontal position $x$} & \multicolumn{6}{c}{horizontal position $y$} & \multicolumn{6}{c}{horizontal position $z$} \\
					\hline
					$\dt$ & $\varepsilon_{L_1}$ & $\mathcal{O}(L_1)$ & $\varepsilon_{L_2}$ & $\mathcal{O}(L_2)$ & $\varepsilon_{L_\infty}$ & $\mathcal{O}(L_\infty)$ & $\varepsilon_{L_1}$ & $\mathcal{O}(L_1)$ & $\varepsilon_{L_2}$ & $\mathcal{O}(L_2)$ & $\varepsilon_{L_\infty}$ & $\mathcal{O}(L_\infty)$ &$\varepsilon_{L_1}$ & $\mathcal{O}(L_1)$ & $\varepsilon_{L_2}$ & $\mathcal{O}(L_2)$ & $\varepsilon_{L_\infty}$ & $\mathcal{O}(L_\infty)$ \\
					\hline
					1.00e-2 & 3.90e-6 & -     & 3.91e-6 & -     & 1.71e-5 & -    & 1.92e-6 & -     & 1.98e-6 & -     & 9.03e-6 & -    & 1.35e-5 & -     & 1.29e-5 & -	   & 4.73e-5 & -    \\
					5.00e-3 & 5.63e-8 & 6.12  & 4.93e-8 & 6.31  & 1.57e-7 & 6.77 & 2.99e-8 & 6.01  & 2.99e-8 & 6.05  & 1.70e-7 & 5.73 & 2.24e-7 & 5.91  & 2.20e-7 & 5.87  & 1.19e-6 & 5.31 \\
					2.50e-3 & 8.54e-10 & 6.04  & 7.33e-10 & 6.07  & 1.27e-9 & 6.94 & 4.58e-10 & 6.03  & 4.32e-10 & 6.11  & 2.77e-9 & 5.94 & 3.47e-9 & 6.02  & 3.26e-9 & 6.07  & 2.07e-8 & 5.85 \\
					1.25e-3 & 1.32e-11 & 6.01  & 1.14e-11 & 6.01  & 1.49e-11 & 6.42 & 7.05e-12 & 6.02  & 6.39e-12 & 6.08  & 4.39e-11 & 5.98 & 5.36e-11 & 6.02  & 4.85e-11 & 6.07  & 3.31e-10 & 5.97 \\
					\hline
					\multicolumn{7}{c}{}  \\
					
				\end{tabular}
				\label{tab:convTable}
			\end{footnotesize}
	\end{table}
\end{center}

%

\subsection{Comparison test against the SPT linear linking approach}
\label{test2}

This test is set up aiming at comparing the SHOT-R method against the SPT linear linking technique, demonstrating an improvement over state-of-the-art linear tracking approaches. In particular, in a time computational domain $T=[t_0;t_f]=[0;2]$, we prescribe the particle velocity given by

\begin{equation}
		v_x(t)= 1 + \tanh(5(t-1)),
		\qquad v_y(t)=2\pi \, \cos(2\pi \,t),
		\label{eqn:test2-ICv}
\end{equation}

from which it is possible to evaluate an analytical expression for the particle trajectory by integrating \eqref{eqn:test2-ICv} over time $t$, thus
\begin{eqnarray}
		x(t)&=& \int_t v_x(t) \, dt = 2t - \frac{\log(\tanh(5(t - 1)) + 1)}{5}, \\	
		y(t)&=&\int_t v_y(t) \, dt=\sin(2\pi \, t).
		\label{eqn:test2-ICx}
\end{eqnarray}

\begin{figure}
	\centering
	\includegraphics[width=1\linewidth]{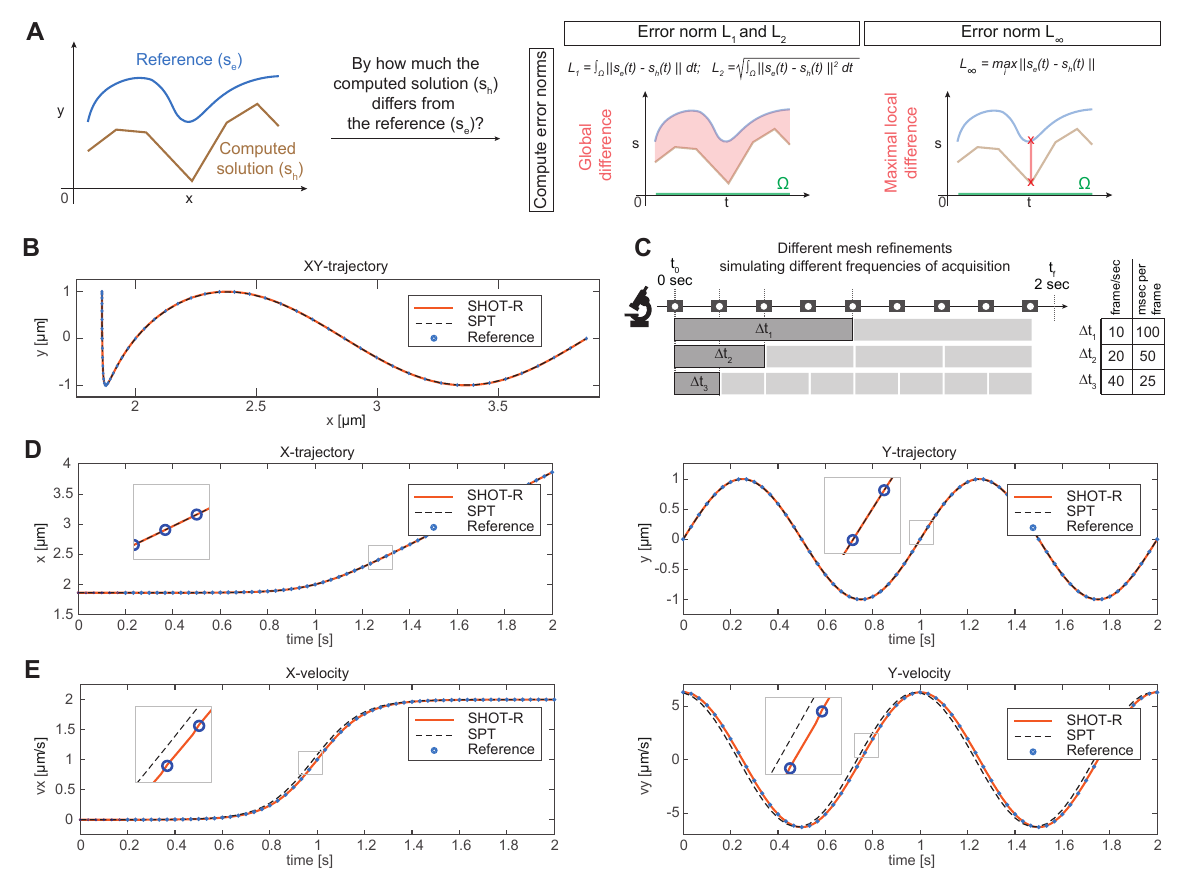}
	\caption{{\bf Comparison test against SPT linear linking.} {\bf (A)} Schematic representation of error norms calculation Eq. \eqref{eqn:Lnorms}. {\bf (B,C)} XY-reference trajectory mimicking particle movement (B) and different mesh refinements simulating timelapse experiments at different acquisition frequency (C). Parameters used in the new mathematical test: Eq.\eqref{eqn:test2-ICv} and Eq.\eqref{eqn:test2-ICx} x- and y- trajectory and velocity in a time computational domain $T=[t_0;t_f]=[0;2]$. {\bf (D,E)} Qualitative comparison between SPT and SHOT-R in trajectory reconstruction (D) and in computing velocity (E) on given equations \eqref{eqn:test2-ICv} and \eqref{eqn:test2-ICx}. (D,E) $\dt_2$ mesh refinement. Abbreviations: $P1$,$P2$,$P3$, first, second, third degree polynomial; $\dt$, time intervals (mesh refinement); SHOT-R, Spatiotemporal High-Order Trajectory Reconstruction; SPT, single particle tracking.}
	\label{fig:Fig2}
\end{figure}

\begin{figure}[!htbp] 
	\centering
	
	\includegraphics[width=1\linewidth]{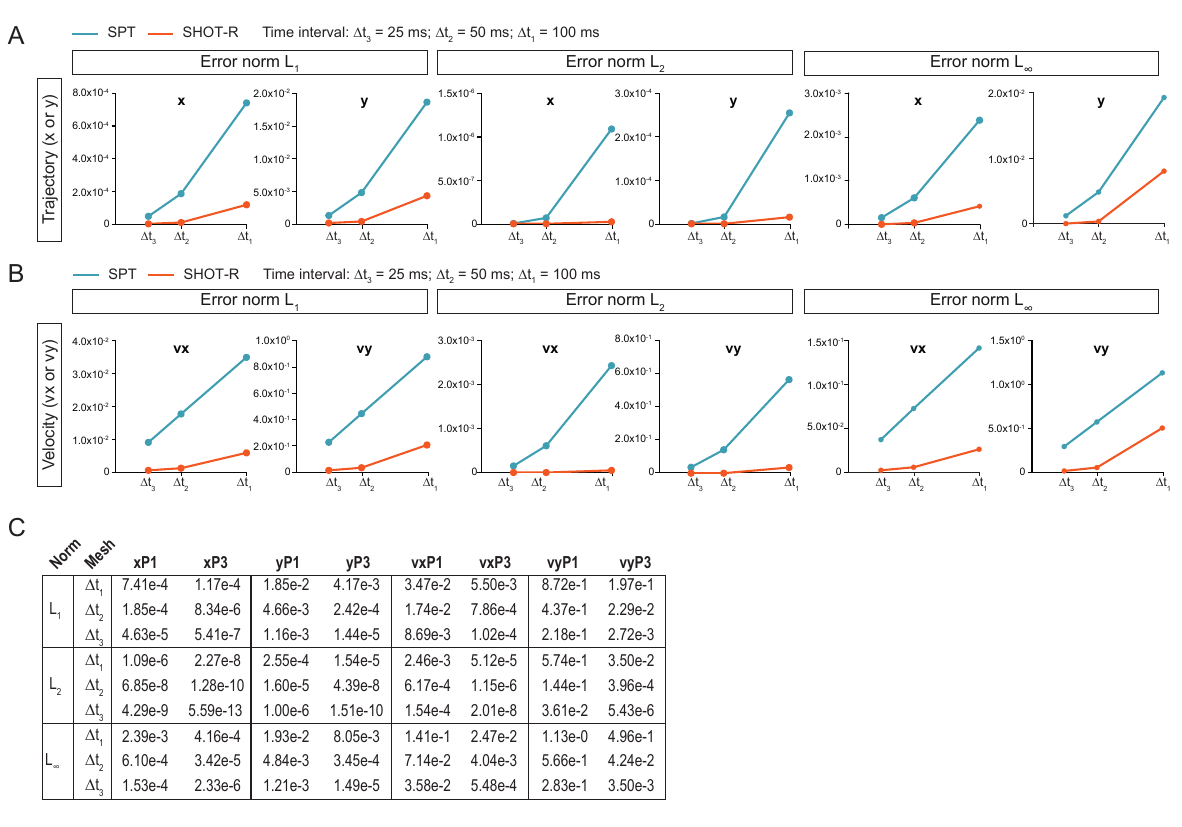}
	
	\caption{{\bf Quantitative comparison test against SPT linear linking.}  {\bf (A-C)} Quantitative comparison between SPT and SHOT-R in trajectory reconstruction (A,C) and in computing velocity (B,C) on given equations Eq.\eqref{eqn:test2-ICv} and Eq.\eqref{eqn:test2-ICx}. Abbreviations: $\dt$, time intervals (mesh refinement); SPT, single particle tracking.}
	\label{fig:Fig3}
\end{figure}

As such, we describe through known functions \eqref{eqn:test2-ICv}-\eqref{eqn:test2-ICx} the shape of an hypothetical complex $xy$-trajectory: along the $x$-axis the particle stalls, then accelerates and eventually moves at a constant velocity; while along $y$-axis it oscillates regularly (Fig. \ref{fig:Fig2}B-E). Three different mesh refinements are considered ($\dt_1$, $\dt_2$ and $\dt_3$), for a total number of $N_1=21$, $N_2=41$ and $N_3=81$ equidistant points of coordinates $\x_k=(x_k,y_k)$, simulating increasing frequencies of acquisition (increased input points $\x_k$) and thus increasing temporal resolution (Fig. \ref{fig:Fig3}C). The particle trajectory is reconstructed and the velocity computed with both the linear linking (SPT) and SHOT-R. 

While with an intermediate mesh density ($\dt_2$) the position is overall nicely reconstructed by both SPT and SHOT-R, the velocity vector is more accurate using a higher order reconstruction (Fig. \ref{fig:Fig2}B-E, \ref{fig:Fig3}A-C). Noteworthy, in the reconstruction of the $v_x$ component, the SPT solution is not capable of reproducing the curve accurately, specifically where the velocity drastically increases (Fig. \ref{fig:Fig2}E). As measured in $L_1$, $L_2$ and $L_\infty$ norms according to \eqref{eqn:Lnorms} (Fig. \ref{fig:Fig2}A), the SHOT-R reconstruction systematically achieves better results with errors of about two orders of magnitude smaller with respect to the linear reconstruction SPT (Fig. \ref{fig:Fig3}A-C). Moreover, the SHOT-R method results to be consistent and convergent, with decreasing errors by increasing the computational mesh resolution (i.e. the more known points $\x_k$ are given as input, the lower are the errors generated by the reconstruction procedure).

{\elo In conclusion, the convergence analysis and the mathematical test validate the accuracy and the robustness of SHOT-R, and show its increased accuracy compared to SPT in computing  velocity.}

\subsection{Velocity reconstruction and backward trajectory integration on biological datasets}
\label{test4}

\begin{figure}[!htbp]
	\centering
	\includegraphics[width=0.7\linewidth]{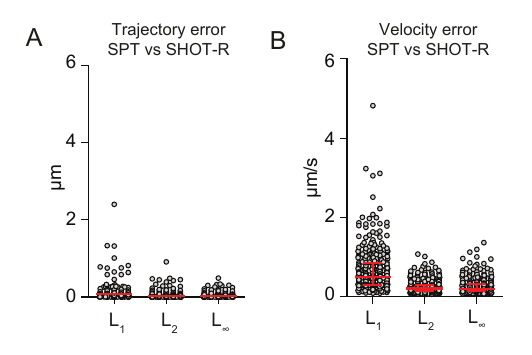}
	\caption{{\bf Comparison of SPT and SHOT-R trajectories and computed velocities on a biological dataset.} {\bf (A,B)} Difference between SPT and SHOT-R reconstructed trajectory (A) or computed velocity (B) measured as error norms ($L_1$, $L_2$, $L_\infty$) with \eqref{eqn:Lnorms}. Abbreviations: SHOT-R, Space-time High-Order Trajectory Reconstruction; SPT, single particle tracking. {\elo Data information: each data point corresponds to a particle. Total number of analyzed particles: 372.}}
	\label{fig:Fig4}
\end{figure}

\begin{figure}[!htbp] 
	\centering
	
	\includegraphics[width=0.8\linewidth]{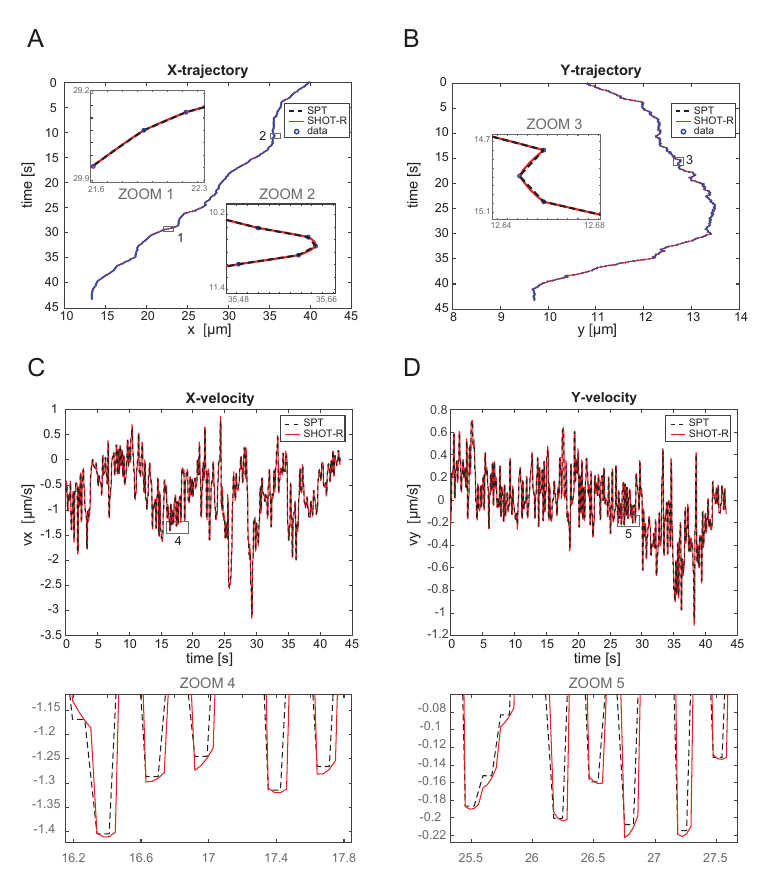}
	
	\caption{{\bf Velocity reconstruction on a biological dataset.} {\bf (A-D)} Representative comparison of SPT and SHOT-R in $xy$ trajectory reconstruction (A,B) and in computing velocity (C,D).}
	
	\label{fig:Fig5}
\end{figure}

As expected, SHOT-R exhibits higher resolution in the reconstruction of the velocity on a mathematical test case (Fig. \ref{fig:Fig2}E, \ref{fig:Fig3}A-C), thus finally we investigate whether our numerical model applied to a biological dataset will produce similar results. In particular, the novel method is tested on a recent study of pre-miR-181a-1 trafficking along axons \cite{corradi2020}. Previously, pre-miR-181a-1 kinematics was studied by kymograph analysis \cite{corradi2020}, while here we retrieve particle coordinates with TrackMate \cite{tinevez2017trackmate} and apply SHOT-R method. SHOT-R trajectory reconstruction is compared with the already existing SPT followed by linear linking methods (Fig. \ref{fig:Fig4}A,B). As observed for the comparison test (Fig. \ref{fig:Fig2}E, \ref{fig:Fig3}A-C), biological particle trajectories do not change depending on the polynomial reconstruction (Fig. \ref{fig:Fig4}A, Fig. \ref{fig:Fig5}A,B). Indeed, the mean values of norm distribution were close to zero (i.e. $L_1$, 0.08 $\mu m$; $L_2$ and $L_\infty$, 0.03 $\mu m$) (Fig. \ref{fig:Fig4}A), suggesting a perfect overlap of the SPT and SHOT-R reconstructed trajectories (Fig. \ref{fig:Fig5}A,B). 
On the contrary, the computed velocities differ with mean errors greater than zero ($L_1$, 0.65 $\mu m/s$; $L_2$, 0.24 $\mu m/s$; $L_\infty$, 0.26 $\mu m/s$) (Fig. \ref{fig:Fig4}B, Fig. \ref{fig:Fig5}C,D). 
Axonal particles moving faster than 0.5 $\mu m/s$ are considered fast-moving puncta \cite{corradi2020,maday2014,leung2018}, hence such an average difference between the two approaches might interfere with the final data interpretation. 

\begin{figure}[!htbp] 
	\centering
	\includegraphics[width=0.7\linewidth]{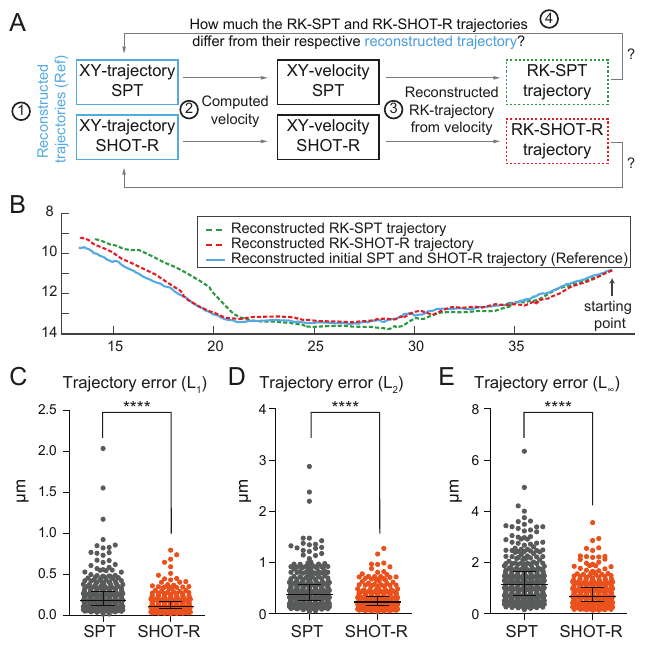}
	\caption{{\bf Backward trajectory integration on a biological dataset.} {\bf (A)} Schematic of the pipeline to test SPT and SHOT-R accuracy in calculating velocity. Color code matches panel B. {\bf (B)} Representative comparison of RK-SPT and RK-SHOT-R trajectory obtained by backward integration of the reconstructed velocity of a retrograde moving particle. {\bf (C-E)} Difference between RK-SPT and RK-SHOT-R measured as error norms $L_1$ (C), $L_2$ (D), $L_\infty$ (E) of the reconstructed trajectory based on velocity trace. Abbreviations: SHOT-R, Spatiotemporal High-Order Trajectory Reconstruction; SPT, single particle tracking; RK, Runge-Kutta method for backward time integration. Data information: **** P $<$ 0.0001. Values are median with interquartile range. Data are not normally distributed (Shapiro-Wilk test), two-tailed Wilcoxon matched pairs test. Each data point corresponds to a particle. Total number of analyzed particles: 372 (C-E).}
	\label{fig:Fig6}
\end{figure}

Investigating which method better approximates the real particle trajectory, and which one is more reliable in studying the particle kinematics remains a major challenge. Indeed, any experimental dataset is always a list of space-time coordinates without a continuous information of particle position during time. As such, it does not exist a reference solution, a known real kinematic behavior of a particle, and thus so far, there has not been a strategy to determine which analysis approach is more accurate over another one on a real experimental dataset. Here, to overcome the problem of the lack of an exact reference solution, we implement a test to check whether our model is significantly more accurate in computing velocities than available strategies (Fig. \ref{fig:Fig6}A). We underline that this approach can be applied to any method for particle kinematic analysis and any kind of experimental dataset.

First, $xy$ trajectories of all particles are reconstructed with first and third polynomial degree (i.e. SPT and SHOT-R trajectories)  based on the $\x_k$ coordinates of all tracked particles (Fig. \ref{fig:Fig6}A, \textcircled{1}). Second, velocity in $xy$ is derived from SPT and SHOT-R trajectories (Fig. \ref{fig:Fig6}A, \textcircled{2}). Third, the particle trajectories are then retrieved from the computed velocity by backward time integration of the ODE \eqref{eqn:ODE} through Runge-Kutta (RK) methods \cite{LeVequeODE} obtaining the "RK-trajectories" (Fig. \ref{fig:Fig6}A \textcircled{3}, Fig. \ref{fig:Fig6}B). Lastly, the reconstructed trajectories based on velocities for SPT and SHOT-R (i.e. RK-SPT and RK-SHOT-R trajectories) are compared to the real trajectory given by the known data points $\x_k$ (Fig. \ref{fig:Fig6}A \textcircled{4}, Fig. \ref{fig:Fig6}B). 
The difference between the reference trajectory and the RK-SHOT-R trajectory, measured as error norm, is significantly reduced compared to the difference between the reference trajectory and the RK-SPT (Fig. \ref{fig:Fig6}C-E), both considering globally the curve (Fig. \ref{fig:Fig6}C,D), as well as the maximal difference in position obtained for each trajectory (Fig. \ref{fig:Fig6}E). This means that the velocities computed by SHOT-R are more accurate than those obtained by SPT, indeed the trajectory retrieved by backward time integration of SHOT-R velocities are more accurate than those retrieved with SPT (Fig. \ref{fig:Fig6}C,E). The Runge-Kutta schemes used for carrying out the time stepping are reported in \ref{app.RK}.

{\elo In summary, we have shown that, likewise for the mathematical test, SHOT-R method drastically increases the accuracy in velocity reconstruction on a real biological dataset.}

\section{Conclusions} \label{sec.concl}
In live-imaging, a series of frames capturing particle position over time is acquired. Retrieving trajectories from these images means approximating the entire particle path. Any linear linking algorithm describes a particle path by joining the known spatial coordinates with straight lines. Here, we reconstruct particle trajectories using polynomials of higher accuracy than second order (i.e. simple straight lines) to better simulate the original trajectory under study.

The reconstruction process is based on the construction of local high order polynomials within each time interval, that are expressed in terms of an expansion with modal Taylor basis functions. To control the spurious oscillations that may arise from the interpolation of discontinuous datasets, we rely on a CWENO limiting strategy. Differently from spline interpolations, that require the solution of a global system on the entire computational grid, we solve local overdetermined linear systems with constraints that ensure the resulting polynomial to be continuous at cell boundaries. 

We face two different challenges: i) the reconstruction of 3D curvilinear trajectories in time (i.e. solving a 4D computational reconstruction), ii) the demonstration of more accurate results using SHOT-R over SPT and linear linking considering the lack of a reference solution in real live-imaging experiment. In our method, we solve these two problems, making them two strengths of SHOT-R itself.
First, to overcome the computational challenge of solving a 4D problem, SHOT-R reconstructs particle trajectories starting from the raw space-time coordinates and splitting the trajectories along each spatial direction ($x$, $y$, $z$), then it reassembles the 4D information by gathering the reconstruction along each axis in 1D without sacrificing the spatial or temporal resolution. The second challenge is approached by measuring the difference between the initial reconstructed trajectory and the one retrieved from the velocities by backward time integration (Fig. \ref{fig:Fig6}). The more accurate the velocity the smaller the difference between the two trajectories. We consider this as an important test to quantify how much our method gains in accuracy compared to SPT in any tested dataset. Indeed, when the particles under study do not move or move at a constant velocity, SHOT-R is not expected to bring more information than what the classical SPT and linear linking do (Fig. \ref{fig:Fig2}). On the other hand, when the moving object (e.g. molecules, virus, vesicles, organelles, or whole cells) shows changes in direction and/or velocity, SHOT-R improves the accuracy in the analysis. 

Beside checking the accuracy of SHOT-R over SPT current strategies, it is crucial to understand which are the new information that can be extracted using the novel algorithm. Specifically, SHOT-R is based on high order and piecewise continuous trajectory reconstructions, hence broadening the range of information that can be unraveled on particle kinematics. Indeed, a higher order method improves the accuracy in computing particle velocities and {\elo would allow} to investigate novel kinematics features (i.e. acceleration). Moreover, {\elo the SHOT-R} continuum trajectory reconstruction {\elo would enable} i) to assess instantaneous velocity and acceleration at any space-time point, ii) to define phenomena based on time independently of the frequency in acquisition (e.g. pausing defined for a temporal threshold in seconds and not in frames), iii) to study co-trafficking events by automatically identify particles moving together, defining a common reference trajectory thus enabling to uncover the kinematics of co-trafficking. 

Our numerical scheme is extremely versatile: it {\elo could} be applied to any live-imaging experiment where the space-time coordinates can be retrieved. In order to address different biological questions, SHOT-R {\elo could provide} an overview on particle dynamics based on the entire trajectory or portions of trajectory, from average to instantaneous behavior, accounting for complex bidirectional motions interspersed with pauses. Moreover, it {\elo would enable} to study particle directionality along different reference vectors depending on the biological questions, and to investigate the kinematics of co-trafficked particles.

Future research will concern the inclusion of uncertainty quantification that is linked to the acquisition process. To allow for a non-intrusive analysis, one might employ collocation methods to account for the uncertainty in the spatial coordinates, as recently proposed in the context of epidemic data \cite{Covid2}.

Overall we envision that SHOT-R could complement and push forward the state-of-the-art strategies for particle kinematics analysis based on live-imaging experiments.

\section*{Declaration of competing interest}
The authors declare that the research was conducted in the absence
of any commercial or financial relationships that could be construed as a potential conflict of interest.

\section*{Data availability}
Data will be made available upon reasonable request.

\section*{Code availability}
The code is implemented in \textit{Matlab} programming language, version \textit{R2020a}. The source code is available from the corresponding authors upon request, and it will be available upon publication at \textit{https://gitlab.com/Boscheri/shot-r}.

\section*{Author contributions}
{E.C. and W.B. designed the research, implemented and validated the method; W.B and M.T. implemented the numerical software; E.C. performed formal analysis and data visualization; E.C., M-L.B. and W.B. discussed concepts and wrote the manuscript; M.-L.B. helped advising the project and provided the biological dataset; W.B., M.-L.B. provided funding.}

\section*{Acknowledgments}

M.-L.B. acknowledges financial support from the grants: Marie Curie Career Integration (618969 GUIDANCE-miR), The G. Armenise-Harvard Foundation, MIUR SIR (RBSI144NZ4) and MIUR PRIN 2017 (2017A9MK4R). \\ W.B. received financial support by Fondazione Cariplo and Fondazione CDP (Italy) under the project No. 2022-1895 and by the Italian Ministry of University and Research (MUR) with the PRIN Project 2022 No. 2022N9BM3N.

\appendix

\section{Particle coordinates} \label{app.bio_data}
In our biological applications, pre-miRNAs particle $xy$ coordinates are extracted from published movies \cite{corradi2020} with TrackMate v5.2.0 \cite{tinevez2017trackmate}. Raw movies are cropped the axonal portion to analyze. The leading tip of the axon, the growth cone, is kept on the right side of each crop, which is not included in the analysis since the high density in particle signal does not allow to recognize single pre-miRNA puncta. Occasionally, background correction with rolling ball pixel size 10 is applied to increase signal to noise ratio in ImageJ. Difference of Gaussian particle detection (DoG) is used as detector in TrackMate. The estimated diameter is set to 0.6 $\mu m$ and the threshold in the TrackMate detection step is fixed from 10 to 400 depending on the signal to noise ratio of each movie. Median filter, reducing the generation of spurious spots in case of salt and pepper noise, and sub-pixel localization options are selected. The tracking settings are: link max distance "1"; gap closing max distance "2" and gap closing max frame gap "4" (for movies with single channel acquisition); link max distance "5"; gap closing max distance "2" and gap closing max frame gap "1" (for movies with both green and red channel).
Splitting or merging events are not expected, hence Simple Linear Assignment Problem (LAP) Tracker is used. No filters are applied to the track step. The automatically tracked lines are manually checked and minor adjustments applied (i.e. remove signal coming from other axons present in the same field; link tracks belonging to the same particle, only when the particle is visible but not detected due to the applied threshold). Prior to analysis a filter on number of consecutive frames per track is applied in order to clear spurious links done by the automatic tracking. For movies on pre-miRNA endogenous trafficking the minimal number of frames is set at 20, for co-trafficking movies at 5.

\section{Computing error norms} \label{app.err_norms}
To quantitatively appreciate the accuracy of the results, the errors between the reconstructed and a reference solution for both position and velocity are measured in $L_1$, $L_2$ and $L_\infty$ norms, that for a generic spatial direction $s$ are
\begin{eqnarray}
	{L_1} &=& \int_{\Omega}\left\|s_{e}(t)-s_{h}(t)\right\| dt, \nonumber \\ 
	{L_2} &=& \sqrt{\int_{\Omega}\left\|s_{e}(t)-s_{h}(t)\right\|^{2} dt}, \nonumber \\
	{L_\infty} &=& \max_{\Omega} \left\|s_{e}(t)-s_{h}(t)\right\|. \label{eqn:Lnorms}
\end{eqnarray}
The reference solution $s_{e}$ depends on the specific application. In particular, it is given by \eqref{eqn:test2-ICx} and \eqref{eqn:test2-ICv} for position and velocity (Fig. \ref{fig:Fig3},\ref{fig:Fig4}) or the $P3$ polynomial (Fig. \ref{fig:Fig6}). Then, $s_h(t)$ denotes the numerical solution computed with SHOT-R (or the linear reconstruction in Fig. \ref{fig:Fig3},\ref{fig:Fig4}). The integrals appearing in \eqref{eqn:Lnorms} are evaluated using Gaussian quadrature formulae of suitable accuracy (see \cite{stroud}).

\section{Backward integration of particle trajectory}
\label{app.RK}
The results shown in Fig. \ref{fig:Fig6}C-E have been obtained by integrating backward in time the trajectory of the particle using the SHOT-R strategy. From the mathematical viewpoint, this corresponds to the solution of the ODE
\begin{equation}
	\frac{d \x}{d t} = -\vv(\x,t), \qquad \x_0 = \x_{N_K}, \qquad t \in [0;t_{N_K}-t_1],
	\label{eqn.ODEback}
\end{equation}
with $\mathcal{D}=2$ and therefore $\x=(x,y)$ and $\vv=(v_x,v_y)$. The initial condition $\x_0$ of the ODE is given by the $xy$ position of the particle at the end of the trajectory, which coincides with the last input point $N_K$. We proceed as follows. Firstly, SHOT-R method is employed for reconstructing the particle trajectory, i.e. for obtaining $p_i^N(\x,t)$ for all cells $i$. Secondly, the velocity vector $\vv_i$ is computed with \eqref{eqn:evalV}, hence obtaining a high order velocity field defined within each computational cell. Finally, \eqref{eqn.ODEback} is solved up to the final time $t_f=t_{N_K}-t_1$ in such a way that the position of the particle at the beginning of the trajectory is retrieved. Obviously, due to numerical approximation, this will never coincide with the input coordinates of the first point $\x_1$. However, we have shown that higher order reconstructions provide a more accurate results, that is the final reconstructed coordinates which arise from the solution of the ODE \eqref{eqn.ODEback} are closer to the original first point of the trajectory taken as input for the trajectory reconstruction. 

In order to solve \eqref{eqn.ODEback} any classical ODE integrator can be used, provided that it satisfies at least the same accuracy property of the polynomial reconstruction. Here, we rely on a very popular numerical technique, namely the class of Runge-Kutta (RK) schemes \cite{LeVequeODE}. For $P1$ reconstruction (SPT) we use a second order RK method, while for $P3$ (SHOT-R) the fourth order version is adopted. Runge-Kutta schemes explicitly write{\tav , for every intermediate step,} as follows:
\begin{itemize}
	\item second order Runge-Kutta (RK2)
	\begin{eqnarray}
		\x_{n+1} &=& \x_n + \frac{d\tau}{2} \left( \mathbf{k}_1 + \mathbf{k}_2 \right), \\
		\mathbf{k}_1 &=& -\vv(\tau_n,\x_n), \nonumber \\
		\mathbf{k}_2 &=& -\vv\left(\x_n + \frac{d\tau}{2}\mathbf{k}_1,\tau_n+\frac{d\tau}{2}\right); \nonumber 
	\end{eqnarray} 
	\item fourth order Runge-Kutta (RK4)
	\begin{eqnarray}
		\x_{n+1} &=& \x_n + \frac{d\tau}{6} \left( \mathbf{k}_1 + 2\mathbf{k}_2 + 2\mathbf{k}_3 + \mathbf{k}_4\right), \\
		\mathbf{k}_1 &=& -\vv(\tau_n,\x_n), \nonumber \\
		\mathbf{k}_2 &=& -\vv\left(\x_n + \frac{d\tau}{2}\mathbf{k}_1,\tau_n+\frac{d\tau}{2}\right), \nonumber \\
		\mathbf{k}_3 &=& -\vv\left(\x_n + \frac{d\tau}{2}\mathbf{k}_2,\tau_n+\frac{d\tau}{2}\right), \nonumber \\
		\mathbf{k}_4 &=& -\vv\left(\x_n + d\tau\mathbf{k}_3,\tau_n+d\tau\right). \nonumber 
	\end{eqnarray} 
\end{itemize}  
\noindent 
The time step $d\tau$ is used to integrate the ODE \eqref{eqn.ODEback} within the time interval $[0;t_{N_K}-t_1]$ and the superscript $n$ denotes the numerical solution at a given time step. $\tau \in [0;t_{N_K}-t_1]$ represents the local time coordinate which advances the solution in time. We choose $d\tau=0.5$ which corresponds to approximately $3.5$ times the time frame of acquisition of the biological dataset used as proof of concept (i.e. 0.144 $sec$). This has been deliberately chosen in order to highlight even more the benefit of a high order reconstruction, because it contains much more information and thus can produce very good results even with large time steps which cover more than three consecutive frames, i.e. more than three adjacent computational cells. If classical linear tracking is used, i.e. $P1$ reconstruction is applied, the results suffer from lack of resolution (Fig. \ref{fig:Fig6}) and the accuracy is drastically reduced.

\newpage

\bibliographystyle{naturemag}
\bibliography{Bibliography_SHOT-R}

\end{document}